\documentclass[leqno,11pt]{amsart}
\usepackage{amsmath}
\usepackage{amssymb,amscd}
\usepackage{amsthm}
\usepackage{latexsym}
\usepackage[myheadings]{fullpage}
\usepackage{color}

\usepackage{enumitem}

\newtheorem{theorem}{Theorem}[section]

\newtheorem{definition}[theorem]{Definition}
\newtheorem{proposition}[theorem]{Proposition}
\newtheorem{corollary}[theorem]{Corollary}
\newtheorem{remark}[theorem]{Remark}
\newtheorem{notation}[theorem]{Notation}
\newtheorem{example}[theorem]{Example}
\newtheorem{question}[theorem]{Question}

\newcommand{\Z}{{\mathbb Z}}
\newcommand{\N}{{\mathbb N}}
\newcommand{\m}{\mathbf{m}}

\def\opn#1#2{\def#1{\operatorname{#2}}} 
\opn\Ap{Ap}
\opn\PF{PF}
\opn\F{F}
\opn\tr{tr}
\opn\res{res}

\address{Dumitru I. Stamate - Faculty of Mathematics and Computer Science, University of Bucharest -  Str. Academiei 14 - 010014 - Bucharest - Romania}
\email{dumitru.stamate@fmi.unibuc.ro}

\address{Francesco Strazzanti - Dipartimento di Scienze Matematiche e Informatiche, Scienze Fisiche e Scienze della Terra - Universit\`a di Messina - Viale Ferdinando Stagno d’Alcontres 31 - 98166 - Messina - Italy.}
\email{francesco.strazzanti@unime.it}

\dedicatory{To the memory of Professor J\"urgen Herzog, whose ideas continue to inspire}

\title{Nearly Gorenstein and almost symmetric properties in shifted numerical semigroups}

\author{Dumitru I. Stamate, Francesco Strazzanti}

\subjclass[2010]{20M25, 20M14, 13H10}

\date{}

\begin{document}

\keywords{Numerical semigroup, shifted family, Frobenius number, pseudo-Frobenius numbers, nearly Gorenstein, almost Gorenstein, reduced type}

\begin{abstract} Given the integers $0<r_1<\dots<r_k$, we consider the shifted family of semigroups $M_n=\langle n, n+r_1,\dots, n+r_k\rangle$, where $n>0$.  For sufficiently large $n$, we prove that if $M_n$ is nearly Gorenstein or almost symmetric, then so is $M_{n+r_k}$. A key ingredient is to relate the pseudo-Frobenius elements of $M_n$ and $M_{n+r_k}$, correcting a wrong claim in the literature. Moreover, we derive explicit formulas for the Frobenius and pseudo-Frobenius numbers of $M_{n+r_k}$.
\end{abstract}

\maketitle

\section*{Introduction}

Let $H\subseteq \N$ be an additive semigroup containing zero and let $\Bbbk[H]=\Bbbk[t^h\mid h\in H]$ denote the associated semigroup ring, where $\Bbbk$ is any field and $t$ is an indeterminate. It has always been of interest to relate algebraic properties  of $\Bbbk[H]$ to those of $H$. With few exceptions, it is generally difficult to find explicit formulas for the invariants of $\Bbbk[H]$ in terms of the generators of $H$, even though general algorithms are available. 

The monoid generated by the nonnegative integers $r_0<r_1<\dots < r_k$ is $\langle r_0, \dots, r_k\rangle=\{\sum_{i=0}^k \lambda_i r_i \mid \lambda_i \in \N, i=0,\dots, k\}$.  One fruitful idea, prompted by a conjecture of Herzog and Srinivasan, was to consider the family of semigroups obtained from it by increasing all generators by the same amount. Set $M_n=\langle n+r_0, n+r_1, \dots, n+r_k\rangle$, where $n\geq 0$. We refer to $M_n$ as a shifted semigroup, and to $\{ M_n \}_{n\geq 0}$ as the shifted family of semigroups obtained from the $r_i$'s, usually omitting the curly brackets when there is no risk of confusion.

Herzog and Srinivasan conjectured the following result on the Betti numbers of these semigroup algebras, later proved by Vu.

\medskip

\noindent{\bf Theorem} (Vu, \cite[Theorem 1.1]{V}) $\beta_i(\Bbbk[M_n])=\beta_i(\Bbbk[M_{n+(r_k-r_0)}])$ for all $i$ and $n\gg 0$.

\medskip

Similar periodicity results hold asymptotically for the Betti numbers of the tangent cones (\cite{HS-cone}) and the homogeneous coordinate ring of the associated projective monomial curves (\cite{V}). More periodic properties have been investigated in \cite{CGHOPPW, JSZ-homog, HO, S-3semi, HHS1, CS, S-CI}.

One consequence of Vu's result is that the complete intersection property and the Gorenstein property  for these algebras eventually occur with period $r_k-r_0$ for $n\gg 0$. The statement about complete intersections had been confirmed independently by Jayanthan and Srinivasan in \cite{JS}.

The main purpose of this paper is to reveal whether similar periodicity occurs for other properties (nearly Gorenstein, almost Gorenstein) and invariants (the residue, the reduced type). When studying asymptotic properties in the shifted family, for simplicity we always assume that $r_0=0$ and that $M_n$ is a numerical semigroup, i.e. $\N\setminus M_n$ is finite, or equivalently that $n$ and $\gcd(r_1, \dots, r_k)$ are coprime.

For a numerical semigroup $H$, its pseudo-Frobenius numbers are the elements of the set $\PF(H)=\{x\in \Z\setminus H \mid x+h\in H \text{ for all } 0\neq h\in H\}$. The Frobenius number of $H$ is $\F(H)=\max \PF(H)= \max \mathbb{Z} \setminus H$. 

The set $\PF(M_n)$ for large $n$ has been investigated by O'Neill and Pelayo in \cite{OP}. They assert that a certain map $\varphi_n: \PF(M_n) \to \PF(M_{n+r_k})$ is a bijection, but we find counterexamples to this statement, see Example \ref{first example}. In Theorem \ref{PF} we give the correct bijection $\varphi_n$. 
This allows us to derive formulas for the Frobenius and pseudo-Frobenius numbers of $M_{n+\lambda r_k}$ for every $\lambda \in \mathbb{N}$ based on those of $M_n$, see Corollary \ref{lambda} and Proposition \ref{Frobenius}.
One by-product is in Proposition \ref{PF ordered}: if $f_1<\dots <f_t$ are the pseudo-Frobenius numbers of $M_n$, then $\varphi_n(f_1)<\dots<\varphi_n(f_t)$ are the pseudo-Frobenius numbers of $M_{n+r_k}$, for $n\gg 0$.

According to \cite{HHS}, a Cohen-Macaulay local ring $(R,\mathfrak{m})$ or a positively graded $\Bbbk$-algebra with maximal graded  ideal $\mathfrak{m}$ and canonical module $\omega_R$ is called a {\em nearly Gorenstein ring} if the trace ideal of the canonical module satisfies $\tr(\omega_R)\supseteq \mathfrak{m}$. It is known that the ring $R$ is Gorenstein if and only if $\omega_R\cong R$, or equivalently $\tr(\omega_R)=R$. We refer to \cite{BH} and \cite{HHS} for the unexplained terminology and background. If $\Bbbk[H]$ is a nearly Gorenstein ring, we say that the monoid  $H$ is nearly Gorenstein.
With notation as above, we prove the following:

\medskip
\noindent {\bf Theorem \ref{NearlyGorenstein}.} For $n\gg 0$, if $M_n$ is nearly Gorenstein, then $M_{n+r_k}$ is nearly Gorenstein, too.
\medskip

This was confirmed when $k=2$ in \cite[Corollary 4.5]{HHS1}, relying  on a good understanding of the asymptotic trace ideals  when $k$ is small (\cite{S-3semi}).

Moscariello and Strazzanti (\cite{MS}) show that if $H$ is a numerical semigroup generated by $h_1,\dots, h_s$, then $\Bbbk[H]$ is a nearly Gorenstein ring if and only if there exists a vector $(f_1,\dots, f_s)\in \PF(H)^s$ such that $h_i+f_i-f\in H$ for all $f\in \PF(H)$ and all $i$. Such a vector is called a {\em nearly Gorenstein vector} for $H$. In fact, what we prove in Theorem \ref{NearlyGorenstein} is that if $(f_0,\dots, f_k)$ is a nearly Gorenstein vector for $M_n$, then so is $(\varphi_n(f_0),\dots, \varphi_n(f_k))$ for $M_{n+r_k}$.

This more general result has further consequences. The almost Gorenstein property was introduced by Barucci and Fr\"oberg for analytically unramified rings in \cite{BF}, and later extended in \cite{Goto1} and \cite{Goto2} to a wider class of Cohen-Macaulay rings in arbitrary dimension.  Given a numerical semigroup $H$, the semigroup ring  $\Bbbk[H]$ is almost Gorenstein precisely when $H$ is an almost symmetric semigroup (cf. \cite{BF}), which in the terminology introduced above means that $(\F(H),\dots, \F(H))$ is a nearly Gorenstein vector for $H$.
Therefore, for $n \gg 0$, it follows that if the semigroup $M_n$ is almost symmetric, then so is $M_{n+r_k}$ (see Corollary \ref{Almost symmetric}).
Moreover, in Proposition \ref{even-type} we show that in this case both have odd type; this fact extends the observation made by Herzog and Watanabe in \cite{HW} (see also \cite[Remark 4.6]{HHS1}) that, when $k=2$, in the shifted family $M_n$ there are only finitely many almost symmetric numerical semigroups that are not symmetric.

For a numerical semigroup $H$ and $R=\Bbbk[H]$, the canonical trace ideal $\tr(\omega_R)$ is generated by the monomials whose exponents are in the semigroup ideal $\tr(H)\subseteq H$.  The residue of $H$ was defined in \cite{HHS1} as $\res(H)= |H\setminus \tr(H)|$. It can be viewed as a measure of how far $H$ is from being symmetric (i.e. $\Bbbk[H]$ is from being a Gorenstein ring): when $\res(H)=0$, the semigroup $H$ is symmetric, and when $\res(H)\leq 1$, $H$ is nearly Gorenstein. It was natural to ask \cite[Question 4.1]{HHS1} whether $\res(M_n)=\res(M_{n+r_k})$ for $n\gg 0$.  The answer is positive when $k=2$, cf. \cite[Theorem 4.2]{HHS1}. However, we provide a counterexample for $k=3$. We show in Example \ref{res} that in the shifted family  $M_n=\langle n, n+2, n+3, n+7 \rangle$ we have $\res(M_{63+7\lambda})=\lambda+9$ for all $\lambda \geq 0$.
We ask whether the correct expectation is that $\res(M_n)$ be a quasilinear function for $n\gg 0$.

In the last section we show that if $\Bbbk[[M_n]]$ has a canonical reduction for $n \gg 0$ (in the sense of \cite{R}), then so does $\Bbbk[[M_{n+r_k}]]$. We also prove that the reduced type of the complete local rings $\Bbbk[[M_n]]$ and $\Bbbk[[M_{n+r_k}]]$ coincide for $n\gg 0$.

For all our results, we provide explicit lower bounds for $n$ from which the periodicity begins. 

\medskip

\noindent{\bf Acknowledgements.}
This work was supported in part by the grant PID2022-137283NB-C22 funded by 
MICIU/AEI/ 10.13039/501100011033 and by ERDF/EU.
The second author was also supported by the ``National
Group for Algebraic and Geometric Structures, and their Applications'' (GNSAGA - INdAM).
Several computations were carried out with the NumericalSgps package \cite{DGM} of the GAP system \cite{GAP}.

\section{Preliminaries} \label{s-prelim}

A {\it numerical semigroup} $H$ is an additive submonoid of $\mathbb{N}$ such that $\mathbb{N} \setminus H$ is finite. If $a_1, \dots, a_s \in \mathbb{N}$, we define $\langle a_1, \dots, a_s \rangle= \{\sum_{i=1}^s \lambda_i a_i \mid \lambda_i \in \mathbb{N}\}$; this is a numerical semigroup if and only if $\gcd(a_1, \dots, a_s)=1$. In this case, we say that $a_1,\dots, a_s$ are generators of $\langle a_1, \dots, a_s \rangle$, and if it is not possible to remove some of them without changing the semigroup, we say that they are {\it minimal generators}. It is well known that every numerical semigroup $H$ has a unique set of minimal generators, and its cardinality is known as the {\it embedding dimension} of $H$. 
The maximum of $\mathbb{Z} \setminus H$ is called the {\it Frobenius number} of $H$ and is denoted by $\F(H)$; note that it could be negative, as $\F(\N)=-1$.
 An integer $f \in \mathbb{Z}\setminus H$ is said to be a {\it pseudo-Frobenius number} of $H$ if $f+h \in H$ for every $h \in H\setminus \{0\}$. The set of pseudo-Frobenius numbers of $H$ is denoted by $\PF(H)$ and its cardinality is called the {\it type} of $H$; of course, $\F(H) \in \PF(H)$, and so the type of $H$ is at least $1$. 

Let $H$ be an additive submonoid of $\N$ (not necessarily a numerical semigroup). Let $d$ be the greatest common divisor of elements in $H$.  Then $H$ is a numerical semigroup precisely when $d=1$. If $d>1$, the Frobenius number of $H$ is defined as $\F(H)= d\F(H')$, where $H'=\{a/d \mid a\in H\}$ is a numerical semigroup. Note that
 $\F(H)=\max\{s \in \Z\setminus H \mid s \equiv 0 \mod d\}$, which also covers the situation when $H$ has a unique minimal generator, i.e. $H'=\N$.

The Ap\'ery sets are an important tool in semigroup theory. Let $H$ be an additive submonoid of $\N$ (not necessarily a numerical semigroup).  For any nonzero element $n\in H$ the  Ap\'ery set of $H$ with respect to  $n$ is defined as $\Ap(H,n)=\{h \in H \mid h-n \notin H\}$.
The Ap\'ery set of a numerical semigroup contains much useful information. For example, the largest element of $\Ap(H,n)$ is $\F(H)+n$.

The monoid $H$ naturally defines a partial ordering on the integers by letting $x\leq_H y$ whenever $y=x+h$ for some $h\in H$.  With respect to this ordering, the maximal elements of $\Ap(H,n)$ are precisely the elements of the form $f+n$ with $f \in \PF(H)$. This provides a method to compute $\PF(H)$.

A factorization of an integer $x$ with respect to the integers $a_1, \dots, a_s$ is an expression of the form $x=\lambda_1 a_1 + \dots + \lambda_s a_s$ with $\lambda_1, \dots, \lambda_s \in \mathbb{N}$; its length is $\lambda_1+\dots+\lambda_s$.  By a factorization of $x$ in the numerical semigroup $H$, we mean a factorization of $x$ with respect to the minimal generators of $H$. For further details on numerical semigroups and unexplained terminology, see \cite{ADG,RG}.

\medskip


\begin{notation}\label{nota}\rm
The following notation will be used throughout the paper.
\begin{itemize}[leftmargin=*]
\item $k\geq 2$, and $0=r_0<r_1<\dots < r_k$ are integers with $d=\gcd(r_1,\dots, r_k)$. 
\item $S=\langle r_1,\dots, r_k\rangle$. 
\item $N_0=\max \{r_k^2, r_k^2+\F(S)r_k\}$.
\item For $n$ a positive integer with $\gcd(n, d)=1$, let $M_n=\langle n,n+r_1,n+r_2, \dots, n+r_k \rangle$. 
\item For $i\in S$, we denote by $\m(i)$ the minimum length of a factorization of $i$  with respect to $r_1, \dots, r_k$. 
\end{itemize}
\end{notation}

\medskip

Note that the monoid $S=\langle r_1,\dots, r_k\rangle$ could have less  than $k$ minimal generators (even only one). 
The condition $\gcd(n,d)=1$ is precisely what is needed for the monoid $M_n$ to be a numerical semigroup. When $n> r_k$, the sum of any two of the listed generators for $M_n$ is larger than $n+r_k$, hence $n, n+r_1,\dots, n+r_k$ generate $M_n$ minimally.

  From the definition, for $i\in S$ it follows that $\m(i)r_1 \leq i \leq \m(i)r_k$.
If $r_1, \dots, r_k$ are the minimal generators of $S$, then $\m(i)$ coincides with the minimum length of a factorization of $i$ in $S$. However, in general, these two values may differ. For example, if $M_n=\langle n, n+3, n+4, n+6\rangle$ then $S=\langle 3,4,6 \rangle=\langle 3,4 \rangle$.  Clearly, $12=2\cdot 6$ yields $\m(12)=2$, while the minimum length of a factorization of $12$ in $S$ equals $3$.

In \cite{OP}, $\m(i)$ is set to denote the minimum factorization length in $S$. However, what is actually used in the proofs is the minimum length of factorizations with respect to $r_1, \dots, r_k$. In this paper, we will make use of some results of O'Neill and Pelayo from \cite{OP}, which remain valid under our definition of $\m(i)$.  We will frequently use their following theorem without explicitly recalling it.

\begin{theorem}{\rm (O'Neill, Pelayo \cite[Theorem 3.3]{OP})} \label{Apéry}
If $n > r_k^2$ and $dn \in S$, then 
\[
\Ap(M_n,n)=\{i+\m(i) n \mid i \in \Ap(S,dn)\}.
\]
Moreover, for each $i \in \Ap(S,dn)$, all the factorizations of $i+\m(i) n$ in $M_n$ have length $\m(i)$.
\end{theorem}

In the following remark, we collect some properties that will be used frequently throughout the paper.

\begin{remark} \rm \label{OP}
(1) By \cite[Theorem 4.3]{BOP}, it follows that $\m(i+r_k)=\m(i)+1$ if $i\in S$ with $i>r_{k-1} r_k$. For instance, this happens when $i >r_k^2-r_k$, or when $i>r_k^2+\F(S)r_k$, since $r_k^2-r_k \geq r_k^2-d r_k=(r_k-d)r_k \geq r_{k-1} r_k$. Notice that in \cite[Theorem 4.3]{BOP} the authors assume that $\{r_1, r_2, \dots, r_k\}$ are the minimal generators of a numerical semigroup, but their proof works also in our context.

(2) As noticed in \cite[Proposition 3.4]{OP}, if $dn>\F(S)$, then $\Ap(S,dn)=\{i_0, \dots, i_{n-1}\}$, where
\[
i_j=
\begin{cases}
dj & \text{ if \ } dj \in S, \\
dj+dn & \text{ if \ } dj \notin S. \\
\end{cases}
\]
Indeed, it is easy to see that $i_j \in \Ap(S,dn)$ and that they are distinct modulo $dn$, therefore it is enough to note that $|\Ap(S,dn)|=n$. 

(3) By the definition of the Apéry set, it is clear that $i\leq  dn+\F(S)$ for every $i \in \Ap(S,dn)$.

(4) Note that $r_k^2 > \F(S)+2r_k$. To see this, set $m=\gcd(r_{k-1},r_k)$ and $T=\langle r_{k-1}/m,r_k/m \rangle$. By the definition of $\F(S)$, it follows that $m\F(T) \geq \F(S)$, whereas the well-known Sylvester formula \cite[Proposition 2.13]{RG} gives $\F(T)=r_{k-1}r_k/m^2-r_{k-1}/m-r_k/m$. Therefore, $r_k^2\geq r_k(r_{k-1}+1)=r_k+r_k r_{k-1}>r_k+r_k r_{k-1}/m-r_{k-1}-r_k+r_k=r_k+m\F(T)+r_k \geq \F(S)+2r_k$.

(5) If $n>r_k^2$, then by the previous point it follows that $dn>r_k^2>\F(S)$ and hence $dn \in S$. Therefore, in Theorem \ref{Apéry}, the assumption $dn \in S$ can be omitted.
\end{remark}
 
\section{Pseudo-Frobenius numbers}

The pseudo-Frobenius numbers of shifted numerical semigroups were studied in \cite{OP}. In particular, in \cite[Theorem 4.8]{OP} a bijection between $\PF(M_n)$ and $\PF(M_{n+r_k})$ is constructed. This construction contains a mistake, and we start this section by correcting it. Indeed, the correspondence in (\ref{bijection}) of the following theorem is defined in a different way, and this will play a crucial role throughout the paper. See also Example \ref{first example}.

\begin{theorem} \label{PF}
Let $P_n$ denote the set 
\[
P_n=\{i \in \Ap(S,dn) \mid f \equiv i \mod n {\rm \ for \ some \ } f \in \PF(M_n)\}.
\]
For $n > N_0$, the map $\psi_n:P_n \rightarrow P_{n+r_k}$ given by
\begin{equation} \label{bijection}
i \mapsto
\begin{cases}
i &{\text if\ } i < dn-r_k \\
i+d r_k &{\text if\ } i \geq dn-r_k
\end{cases}
\end{equation}
is a bijection. In particular, this induces a bijection $\varphi_n:\PF(M_n) \rightarrow \PF(M_{n+r_k})$.
\end{theorem}

\begin{proof}
We divide the proof into three steps.

\medskip

{\bf Step 1: If $i < dn-r_k$, then $i \in P_n \Longleftrightarrow i \in P_{n+r_k}$.} We may assume that $i \in S$, and in this case $i$ is in both $\Ap(S,dn)$ and $\Ap(S,d(n+r_k))$. It follows that 
\begin{align}
i+\m(i) n &\in \Ap(M_n,n), \label{1.1} \\ 
i+\m(i)(n+r_k) &\in \Ap(M_{n+r_k},n+r_k). \label{1.2}
\end{align}
Since $i <dn-r_k$ and $i \in S$, for every $j=0, \dots, k$ we have $i+r_j<dn$ and then $i+r_j$ is in both $\Ap(S,dn)$ and $\Ap(S,d(n+r_k))$. Thus
\begin{align}
i+r_j+\m(i+r_j) n &\in \Ap(M_n,n), \label{1.1a} \\ 
i+r_j+\m(i+r_j)(n+r_k) &\in \Ap(M_{n+r_k},n+r_k). \label{1.2a}
\end{align}

By (\ref{1.1}), $i \in P_n$ if and only if $i+(\m(i)-1) n \in \PF(M_n)$, which means that $i+(\m(i)-1) n+n+r_j=i+r_j+\m(i)n \in M_n$ for every $j=0, \dots, k$. By (\ref{1.1a}), this is equivalent to $\m(i+r_j) \leq \m(i)$ for every $j=0, \dots, k$. 

Analogously, $i \in P_{n+r_k}$ if and only if $i+(\m(i)-1) (n+r_k)+n+r_k+r_j=i+r_j+\m(i) (n+r_k) \in M_{n+r_k}$ for every $j=0, \dots, k$, and also this condition is equivalent to $\m(i+r_j) \leq \m(i)$ by (\ref{1.2a}).

\medskip

{\bf Step 2: If $i \geq dn-r_k$, then $i \in P_n \Longleftrightarrow i+dr_k \in P_{n+r_k}$.} Clearly, $i\geq dn-r_k\geq r_k^2-2r_k>\F(S)$. Therefore, $i \in \Ap(S,dn)$ if and only if $i+d r_k \in \Ap(S,d(n+r_k))$ because $i+dr_k-d(n+r_k)=i-dn$, so we assume this is the case.

As in Step 1, $i \in \Ap(S,dn)$ and $i+dr_k \in \Ap(S,d(n+r_k))$ imply $i+\m(i)n \in \Ap(M_n,n)$ and $i+dr_k+\m(i+dr_k)n \in \Ap(M_{n+r_k},n+r_k)$. Thus, $i \in P_n$ if and only if $i+(\m(i)-1) \, n \in \PF(M_n)$, and then 
\begin{equation}
i \in P_n \Longleftrightarrow i+r_j+\m(i) n \in M_n \text{ for every } j=0, \dots, k. \label{2.1}
\end{equation}
Moreover, since $i\geq dn-r_k>r_k^2-r_k$, Remark \ref{OP}.(1) implies that $\m(i+dr_k)=\m(i)+d$, and then we also have
\begin{equation} \label{2.2}
i+dr_k \in P_{n+r_k} \Longleftrightarrow i+dr_k+r_j+(\m(i)+d)(n+r_k) \in M_{n+r_k} \text{ for every } j=0, \dots, k.
\end{equation}
We fix $j$ and distinguish two cases.

\smallskip

$\bullet$ \ Let $i+r_j \in \Ap(S,dn)$. Then, $i+r_j+\m(i+r_j)n \in \Ap(M_n,n)$ implies that $i+r_j+\m(i) n \in M_n$ precisely when $\m(i+r_j)\leq \m(i)$. 
Moreover, $i+dr_k+r_j-d(n+r_k)=i+r_j-dn \notin S$ means that $i+dr_k+r_j \in \Ap(S,d(n+r_k))$ and then
\[
i+dr_k+r_j+\m(i+dr_k+r_j)(n+r_k)=i+dr_k+r_j+(\m(i+r_j)+d)(n+r_k)
\]
is in $\Ap(M_{n+r_k},n+r_k)$. Hence, by (\ref{2.2}), also $i+dr_k \in P_{n+r_k}$ is equivalent to $\m(i+r_j) \leq \m(i)$.

\smallskip

$\bullet$ \ Let $i+r_j \notin \Ap(S,dn)$. Then, $i+r_j-dn \in S$. Moreover, by (3) and (4) of Remark \ref{OP},
$i+r_j-dn \leq dn+\F(S)+r_k-dn<r_k^2-r_k<dn$, thus $i+r_j-dn$ is in both $\Ap(S,dn)$ and $\Ap(S,d(n+r_k))$. 
Therefore, $i+r_j-dn+\m(i+r_j-dn)n \in \Ap(M_n,n)$ and 
$i+r_j-dn+\m(i+r_j-dn)(n+r_k) \in \Ap(M_{n+r_k},n+r_k)$. 
We now show that $\m(i+r_j-dn) \leq \m(i)+d$. Bearing in mind Remark \ref{OP}(3), the inequality follows because
\begin{gather*}
\m(i)+d = \frac{\m(i)r_k}{r_k}+d \geq \frac{i}{r_k}+d \geq \frac{dn-r_k}{r_k}+d \geq \frac{n}{r_k} > \frac{\F(S)r_k+r_k^2}{r_k}=\\
=\F(S)+r_k \geq \F(S)+r_j \geq i-dn+r_j \geq \m(i+r_j-dn). 
\end{gather*}
As a consequence, $i+r_j+\m(i)n=i+r_j-dn+(\m(i)+d)n \in M_n$ and 
$i+dr_k+r_j+(\m(i)+d)(n+r_k)=i+r_j-dn+(\m(i)+2d)(n+r_k) \in M_{n+r_k}$.

\smallskip

Bearing in mind (\ref{2.1}) and (\ref{2.2}), we immediately get that $i \in P_n$ if and only if $i+dr_k \in P_{n+r_k}$.

\medskip

{\bf Step 3: If $dn-r_k \leq i < dn+(d-1)r_k$, then $i \notin P_{n+r_k}$.}
If $i\not\equiv 0 \mod  d$, the conclusion is clear. We now assume $d$ divides $i$. 
We have $\F(S)<r_k^2-r_k< dn-r_k \leq i$ and $i < d(n+r_k)$, which imply $i \in \Ap(S, d(n+r_k))$ and $x:=i+\m(i)(n+r_k) \in \Ap(M_{n+r_k}, n+r_k)$. We need to show that $x$ is not maximal in $\Ap(M_{n+r_k}, n+r_k)$. Note that $\m(i+r_k)=\m(i)+1$ because $i\geq dn-r_k >r_k^2-r_k$. Since $\F(S) <i< i+r_k < dn+(d-1)r_k+r_k = d(n+r_k)$, we have $i+r_k \in \Ap(S, d(n+r_k))$ and then
\[
x+(n+2r_k)=i+r_k+(\m(i)+1)(n+r_k)=i+r_k+\m(i+r_k)(n+r_k) \in \Ap(M_{n+r_k}, n+r_k).
\]
This means that $x$ is not maximal in $\Ap(M_{n+r_k}, n+r_k)$ because $n+2r_k \in M_{n+r_k}$. It follows that $x-(n+r_k) \notin \PF(M_{n+r_k})$, which easily implies $i \notin P_{n+r_k}$.

\medskip
Putting together the three steps, we immediately get that $\psi_n$ is a bijection. 
\end{proof}

With notation as above, \cite[Theorem 4.8]{OP} states that the map 
 $\psi_n:P_n \rightarrow P_{n+r_k}$ is given by 
\begin{equation} \label{bijection-wrong}
i \mapsto
\begin{cases}
i &{\text if\ } i \leq dn \\
i+r_k &{\text if\ } i > dn.
\end{cases}
\end{equation}

The next examples show how Theorem \ref{PF} works. Notice that they agree with Theorem \ref{PF}, but are inconsistent with \eqref{bijection-wrong}.

\begin{example} \rm \label{first example}
1) Consider the shifted family $M_n=\langle n, n+2, n+6, n+7\rangle$. Therefore, $r_1=2, r_2=6, r_3=7, d=1$, $S= \langle 2,6,7 \rangle=\langle 2,7 \rangle$, $\F(S)=5$ and $N_0=r_3^2+\F(S)r_3=84$. For $n=88$ we have $M_{88}=\langle 88, 90, 94, 95\rangle$ and we find  with  GAP (\cite{GAP}) that  $\PF(M_{88})=\{281, 1141, 1145, 1237 \}$.  Since
\begin{align*}
281=17+3\cdot  88,  &\text{ where } 17\in \Ap(S, 88)  \text{ and } \m(17)=4, &\ \ \ \ \ \ \ \  \\
1141=85+ 12\cdot 88, & \text{ where } 85\in \Ap(S, 88)  \text{ and } \m(85)=13,  &\ \ \\
1145= 89+12\cdot 88,  &\text{ where } 89\in \Ap(S, 88)  \text{ and } \m(89)=13, &\ \ \\
1237= 93+13\cdot 88,  &\text{ where } 93\in \Ap(S, 88)  \text{ and }\m(93)=14, &\ \  
\end{align*}
we obtain $P_{88}=\{17, 85, 89, 93\}$. Moreover, if we consider $M_{88+r_3}=M_{95}=\langle 95, 97, 101, 102\rangle$,  the correspondence \eqref{bijection} in Theorem \ref{PF} gives  $P_{95}=\{17, 92, 96, 100\}$. This, together with Theorem \ref{Apéry}, indicates that $\PF(M_{95})=\{ 302, 1327, 1331, 1430\}$, as it can also be verified with GAP.

On the other hand, according to \eqref{bijection-wrong} one should have $P_{95}=\{17, {\bf 85}, 96, 100 \}$, which is not correct.

2) For the shifted family $M_n=\langle n,n+8, n+12, n+14 \rangle$ we have $r_1=8, r_2=12, r_3=14, d=2$ and $S=\langle 8,12,14 \rangle$ with $\F(S)=18$ and $N_0=r_3^2+\F(S)r_3=448$. For  $M_{449}=\langle 449,457,461,463 \rangle$ we have  $\PF(M_{449})=\{29171, 29636, 30101 \}$. Notice that in this case, the smallest element in $P_{449}=\{ 884, 900, 916\}$  is in fact $dn-r_3=884$.  If we consider $M_{449+r_3}=M_{463}=\langle 463,471,475,477 \rangle$ we obtain with GAP that $\PF(M_{463})=\{ 31007, 31486, 31965\}$. According to Theorem \ref{PF}, we have:
\begin{align*}
&29171 \equiv 884 \text{ mod } 449 \text{  \ with } 884 \in \Ap(S,898), &\ \ \ \ 31007 \equiv 912 \text{ mod } 463 \text{ \ with } 912 \in \Ap(S,926), \\
&29636 \equiv 900 \text{ mod } 449 \text{ \ with } 900 \in \Ap(S,898),  &\ \ \ \ 31486 \equiv 928 \text{ mod } 463 \text{ \ with } 928 \in \Ap(S,926), \\ 
&30101 \equiv 916 \text{ mod } 449 \text{ \ with } 916 \in \Ap(S,898),  &\ \ \ \ 31965 \equiv 944 \text{ mod } 463 \text{ \ with } 944 \in \Ap(S,926). 
\end{align*}
Hence $P_{463}=\{912, 928, 944\}$, confirming Theorem \ref{PF}, whereas \eqref{bijection-wrong} would predict \{884, 914, 930\}.

Note that $29171$ is also congruent to $435$ modulo $449$ but $435 \notin \Ap(S,898)$ because it is odd and not in $S$, so the corresponding element in $P_{449}$ is $884$. 
\end{example}

\begin{remark} \rm
As also noted in \cite{OP}, Theorem \ref{PF} implies that the type of $M_{n}$ is periodic for sufficiently large $n$. This result also follows from Vu's  \cite[Theorem 1.1]{V}, where it is shown that the same holds for all the Betti numbers of the corresponding semigroup ring $\Bbbk[M_n]$. However, Theorem \ref{PF} provides an elementary proof of this result for the case of the last nonzero Betti number, which equals the type of $M_n$.
\end{remark}

Throughout the paper we write $\varphi_n(f)$ to denote the image of $f\in \PF(M_n)$ via the bijection $\varphi_n$ defined in Theorem \ref{PF}. Moreover, we write $\varphi_n^\lambda(f)$ to denote the image of $f$ via the composition $\varphi_{n+(\lambda-1) r_k} \circ \varphi_{n+(\lambda-2) r_{k}}\circ \dots \circ \varphi_{n+r_k} \circ \varphi_{n}$; hence, $\varphi_n^\lambda(f) \in \PF(M_{n+\lambda r_k})$. Similarly, we define $\psi_n(f)$ and $\psi^{\lambda}_n(f)$.

\begin{definition}
Given $P_n$ as in Theorem \ref{PF}, we define the following two sets:
\[
P_n'=\{i \in P_n \mid i <dn-r_k\}, \ \ \ \ \ \ \ \
P_n''=\{i \in P_n \mid i \geq dn-r_k\}.
\]
\end{definition}

In the next remark, we collect some properties that we will frequently use throughout the paper.

\begin{remark} \rm \label{P''} Assume that $n >N_0=\max\{r_k^2, r_k^2+\F(S)r_k\}$. \\
(1) Since $i\geq dn-r_k$ is equivalent to $i+dr_k\geq d(n+r_k)-r_k$, it immediately follows that $P_{n+r_k}'=P_n'$ and $P_{n+ r_k}''=\{i+dr_k \mid i \in P_n''\}$. More generally, we have $P_{n+\lambda r_k}'=P_n'$ and $P_{n+\lambda r_k}''=\{i+\lambda d r_k \mid i \in P_n''\}$ for every $\lambda \in \mathbb{N}$. \\
(2) If $i\in P_n'$, we claim that $i<d(N_0+r_k)-r_k$. Indeed, let $\lambda$ be the largest integer such that $n'=n-\lambda r_k$ is greater than $N_0$. By (1), $i \in P_{n'}'$, then
\[i<dn'-r_k = d(n'-r_k)+dr_k-r_k \leq d N_0+dr_k-r_k,\]
where the second inequality follows from the maximality of $\lambda$. \\
(3) If $i \in P_n'$, then $i<d(r_k^3-r_k^2)$. Indeed, by (2) and Remark \ref{OP}.(4), when $\F(S)>0$ we immediately get that
$i < d(r_k^2+(\F(S)+1)r_k) \leq d(r_k^2+(r_k^2-2r_k)r_k)=d(r_k^3-r_k^2)$. 

When $\F(S)<0$, then $\F(S)=-d$  and $N_0=r_k^2$. Using (2), it is enough to check that  
\begin{eqnarray}
d(r_k^2+r_k)-r_k &\leq& d(r_k^3-r_k^2), \text{equivalently} \nonumber \\
0&\leq& d(r_k^2-2r_k-1)+1, \text{ or } \nonumber \\
0&\leq& d( (r_k-1)^2-2)+1. \label{easy}
\end{eqnarray}
When $d=1$, since $r_k\geq 2$ we derive \eqref{easy}. When $d>1$, $r_k\geq r_1+d \geq 2d\geq 4$, hence \eqref{easy} holds, as well.
\end{remark}

As a consequence of Theorem \ref{PF}, we obtain a formula for the pseudo-Frobenius numbers of $M_{n+r_k}$.

\begin{corollary} \label{CorPF}
Assume $n >N_0$. If $f \in \PF(M_n)$ and $f \equiv i \mod n$ with $i \in \Ap(S,dn)$, then
\[
\varphi_n(f)=
\begin{cases}
f+(\m(i)-1)r_k &{\text if \ } i< dn-r_k, \\
f+(\m(i)+2d-1)r_k+dn & {\text if \ } i\geq dn-r_k. 
\end{cases}
\]
\end{corollary}

\begin{proof}
Bearing in mind that $f=(i+\m(i)n)-n$, the result immediately follows from Theorem \ref{PF} because $\varphi_n(f)=(i+\m(i)(n+r_k))-(n+r_k)$ if $i<dn-r_k$ and $\varphi_n(f)=(i+dr_k+\m(i+dr_k)(n+r_k))-(n+r_k)$ if $i \geq dn-r_k$, where $\m(i+dr_k)=\m(i)+d$ by Remark \ref{OP}.(1). 
\end{proof}

\begin{example} \label{Counterexample} \rm
Let $M_{200}=\langle 200, 202, 207, 211 \rangle$, so $n=200, r_1=2, r_2=7, r_3=11$, $d=1$ and $S=\langle 2,7,11 \rangle$. The numerical semigroup $M_{200}$ has type $7$ and its pseudo-Frobenius numbers are:
\[
\begin{array}{l@{\hspace{6em}}l}
819=19+5\cdot 200-200 &\text{where } 19 \in \Ap(S,200) \text{ and } \m(19)=5, \\
1012=12+6\cdot 200-200 &\text{where } 12 \in \Ap(S,200) \text{ and } \m(12)=6, \\
3791=191+19\cdot 200-200 &\text{where } 191 \in \Ap(S,200) \text{ and } \m(191)=19, \\
3805=205+19\cdot 200-200 &\text{where } 205 \in \Ap(S,200) \text{ and } \m(205)=19, \\
3995=195+20\cdot 200-200 &\text{where } 195 \in \Ap(S,200) \text{ and } \m(195)=20, \\
3999=199+20\cdot 200-200 &\text{where } 199 \in \Ap(S,200) \text{ and } \m(199)=20, \\
4003=203+20\cdot 200-200 &\text{where } 203 \in \Ap(S,200) \text{ and } \m(203)=20.
\end{array}
\]
Note that $11$ is not a minimal generator of $S$, but to compute $\m(i)$ we have to consider the factorizations of $i$ with respect to $2,7$, and $11$: for instance, $\m(19)=5$ but the minimal length factorization of $19$ in $S$ is $7$. We have $P_{200}'=\{12,19\}$, while $P_{200}''=\{191,195,199,203,205\}$, thus by Corollary \ref{CorPF} the pseudo-Frobenius numbers of $M_{211}$ are:
\begin{gather*}
863=819+(5-1)\cdot 11, \\
1067=1012+(6-1)\cdot 11, \\
4211=3791+(19+1)\cdot 11+200, \\
4225=3805+(19+1)\cdot 11+200, \\
4426=3995+(20+1)\cdot 11+200, \\
4430=3999+(20+1)\cdot 11+200, \\
4434=4003+(20+1)\cdot 11+200.
\end{gather*}
\end{example}

\medskip

As in Corollary \ref{CorPF}, a formula can be obtained also for $M_{n+\lambda r_k}$.

\begin{corollary} \label{lambda}
Assume $n >N_0$. Let $f \in \PF(M_n)$ and $f \equiv i \mod n$ with $i \in \Ap(S,dn)$. Then
\[
\varphi_n^\lambda(f)=
\begin{cases}
f+(\m(i)-1)\lambda r_k &\text{ \ if \ } i< dn-r_k, \\
f+(\m(i)+(\lambda+1)d-1) \lambda r_k+\lambda dn & \text{ \ if \ } i\geq dn-r_k.
\end{cases}
\]
\end{corollary}

\begin{proof}
If $i<dn-r_k$, then $i<d(n+\mu r_k)-r_k$ for every $\mu \in \mathbb{N}$, and thus we get the formula by applying $\lambda$ times the first case in the formula of Corollary \ref{CorPF}.

Assume now $i \geq dn-r_k$. We will prove by induction that $\varphi_n^\lambda(f)=f+(\m(i)+(\lambda+1)d-1) \lambda r_k+\lambda dn$ for every $\lambda \in \mathbb{N}$. Clearly, it holds for $\lambda=0$. Note that $\psi_n(i)=i+dr_k \geq d(n+r_k)-r_k$ and it is easy to see that in general $\psi_n^{\mu}(i)=i+\mu dr_k \geq d(n+\mu r_k)-r_k$ for every positive integer $\mu$. Therefore, by Corollary \ref{CorPF}, induction hypothesis, and Remark \ref{OP}.(1), we get that 
\begin{gather*}
\varphi_n^\lambda(f)=\varphi_n^{\lambda-1}(f)+(\m(i+(\lambda-1)dr_k)+2d-1)r_k+d(n+(\lambda-1)r_k)\\
=f+(\m(i)+\lambda d-1)(\lambda-1)r_k+(\lambda-1) dn+(\m(i)+(\lambda-1)d+2d-1)r_k+dn+d(\lambda-1)r_k \\
=f+(\m(i)+(\lambda+1) d-1)(\lambda-1)r_k
+\lambda dn+(\m(i)+(\lambda+1)d-1)r_k \\
=f+(\m(i)+(\lambda+1)d-1) \lambda r_k+\lambda dn. \qedhere
\end{gather*} 
\end{proof}

From the proof of Corollary \ref{lambda} we also obtain the following.
\begin{corollary}\label{shift-m}
Assume $n>N_0$.  For $i\in P_n$ one has
\[
\m(\psi_n^\lambda (i))=
\begin{cases}
\m(i)  &\text{ \ if \ } i \in P_n', \\
\m(i)+\lambda d & \text{ \ if \ } i\in P_n^{''}.
\end{cases}
\]
\end{corollary}

We conclude this section by deriving a formula for the Frobenius number of $M_{n+\lambda r_k}$. As a preliminary step, we first prove a result that will also be applied in the next section.

\begin{proposition} \label{PF ordered}
Assume $n \geq r_k^4$. Let $f_1<f_2<\dots < f_t$ be the pseudo-Frobenius numbers of $M_n$ and let $f_p \equiv i_p \mod n$ with $i_p \in \Ap(S,dn)$ for $1 \leq p \leq t$. Then, there exists an integer $0\leq s \leq t$ such that $i_p \in P_n'$ when $p\leq s$ and $i_p \in P_n''$ when $p>s$. Moreover, $\varphi_n^\lambda(f_1) < \varphi_n^\lambda(f_2) < \dots < \varphi_n^\lambda(f_t)$ for every $\lambda \in \mathbb{N}$.
\end{proposition}

\begin{proof}
We first note that $n>r_k^3>r_k^2+\F(S)r_k$ by Remark \ref{OP}.(4). Fix $f_a=i_a+\m(i_a)n-n$ and $f_b=i_b+\m(i_b)n-n$. We will prove these three statements:
\begin{enumerate}[leftmargin=*]
\item If $i_a \in P_n'$ and $i_b \in P_n''$, then $\m(i_a)<\m(i_b)$ and $f_a<f_b$.
\item If $i_a,i_b \in P_n'$ and $f_a<f_b$, then $\m(i_a) \leq \m(i_b)$. 
\item If $i_a,i_b \in P_n''$ and $f_a<f_b$, then $\m(i_a) \leq \m(i_b)$.
\end{enumerate}
Then, the thesis will follow immediately from these facts, Remark \ref{P''}.(1), and Corollary \ref{lambda}. 
\begin{enumerate}[leftmargin=*]
\item By Remark \ref{P''}.(3), we have
\[\m(i_a)\leq \frac{i_a}{r_1} \leq i_a < dr_k^3-dr_k^2<\frac{dr_k^4}{r_k}-1\leq\frac{dn-r_k}{r_k} \leq \frac{i_b}{r_k}\leq \m(i_b).\] 
Therefore, since $i_a<i_b$, we immediately obtain
$f_{a}=i_a+\m(i_a)n-n < i_b+\m(i_b)n-n=f_b$. 
\item Using again Remark \ref{P''}.(3), we obtain
$i_b< d(r_k^3-r_k^2) < r_k(r_k^3-r_k^2)<r_k^4 \leq n$.
It follows that
\[
\m(i_b)=\frac{f_b-i_b+n}{n}>\frac{f_a-i_b+n}{n}>\frac{f_a-i_a+n-i_b}{n}=\m(i_a)-\frac{i_b}{n}>\m(i_a)-1
\]
and therefore $\m(i_b) \geq \m(i_a)$.
\item We recall that $dn+\F(S)-i_b \geq 0$ by Remark \ref{OP}.(3), then
\begin{gather*}
\m(i_a)=\frac{f_a-i_a+n}{n}<\frac{f_b-(dn-r_k)+n+(dn+\F(S)-i_b)}{n}=\\
=\frac{f_b-i_b+n+r_k+\F(S)}{n}=\m(i_b)+\frac{r_k+\F(S)}{n}<\m(i_b)+1,
\end{gather*}
which implies that $\m(i_a) \leq \m(i_b)$. \qedhere
\end{enumerate}
\end{proof}

When $n \geq r_k^4$, the previous proposition implies that the Frobenius number of $M_n$ is associated to an element of $P_n''$, provided that $P_n''$ is not empty. In the following result we prove that this always happens.

\begin{proposition} \label{Frobenius}
Assume $n \geq r_k^4$ and let $\F(M_{n}) \equiv i \mod n$ with $i \in \Ap(S,dn)$. Then, $i \in P_n''$ and
\[\F(M_{n+\lambda r_k})= \varphi_n^\lambda(\F(M_n))=
\F(M_n)+(\m(i)+(\lambda+1)d-1) \lambda r_k+\lambda dn\]
for every $\lambda \in \mathbb{N}$.
\end{proposition}

\begin{proof}
In light of Proposition \ref{PF ordered} and Corollary \ref{lambda}, we only need to prove that $i \in P_n''$. Thus, assume by contradiction that $i \in P_n'$.
We first notice that, using Remark \ref{P''}.(3), we have the following inequalities:
\begin{gather*}
r_k^2+i+\m(i)r_k \leq r_k^2+i+\frac{i}{r_1}r_k=r_k^2+i\left(1+\frac{r_k}{r_1}\right)<dr_k^2+d(r_k^3-r_k^2)\left(1+\frac{r_k}{r_1}\right) \\
=d\left(r_k^2+r_k^3-r_k^2+\frac{r_k^4}{r_1}-\frac{r_k^3}{r_1}-r_k^4+r_k^4\right)=d\left(-r_k^3(-1+r_k)+\frac{r_k^3}{r_1}(r_k-1)+r_k^4\right) \\
=d\left(\left(\frac{r_k^3}{r_1}-r_k^3\right)(r_k-1)+r_k^4\right) \leq dr_k^4 \leq dn,
\end{gather*}
where the penultimate inequality holds because $(r_k^3/r_1-r_k^3)\leq 0$ and $(r_k-1)>0$. Therefore, using also Remark \ref{OP}.(4), we get $\F(S) < r_k^2-2r_k< r_k^2+i+\m(i)r_k <dn$ and then $r_k^2+i+\m(i)r_k \in \Ap(S,dn)$. Theorem \ref{Apéry} implies  that $r_k^2+i+\m(i)r_k+\m(r_k^2+i+\m(i)r_k)n \in \Ap(M_n,n)$ but we have 
\begin{gather*}
r_k^2+i+\m(i)r_k+\m(r_k^2+i+\m(i)r_k)n \\
> i+\m(r_k^2+i+ \m(i)r_k)n=  
i+\m(r_k^2+i)n+\m(i)n  \\
>  i+\m(i)n=\F(M_n)+n.
\end{gather*}
This means that $\F(M_n)+n$ is not the largest element in $\Ap(M_n,n)$, which is a contradiction.
\end{proof}

\begin{remark} \rm
In the proof of Proposition \ref{Frobenius}, we have seen that $P_n'' \neq \emptyset$. On the other hand, $P_n'$ can be empty, as is shown by the shifted family in Example \ref{first example}.(2).
\end{remark}

\section{Nearly Gorenstein property}

In this section we establish certain asymptotic properties in the shifted family of numerical semigroups
 $M_n=\langle n, n+r_1,\dots, n+r_k \rangle$. 
 To obtain a bound on $n$ such that these asymptotic properties hold, we introduce an integer $N$. The reason for this technical definition will become clear in the proof of Theorem \ref{NearlyGorenstein}.

\begin{definition} \label{def N}
Consider $M_n$ with $n>N_0=\max \{r_k^2, r_k^2+\F(S)r_k\}$. We define the following non-negative integers:
\begin{gather*}
N_1=\max\{0,r_k+\ell-i, (\m(\ell)-\m(i))r_k-\ell+i \mid i,\ell \in P_n'\}, \\
N_2=\max\{0,(\m(\ell)-\m(i))r_k-\ell+i \mid i,\ell \in P_n''\},\\
N_3=\max\{0,(\m(\ell)r_k+3r_k)/d \mid \ell \in P_n'\}, \\
N=\max\{N_0+r_k+\F(S)/d,N_1,N_2,N_3\}.
\end{gather*}
\end{definition}

\begin{remark} \rm \label{N}
We claim that the integer $N$ is the same for $M_{n+\lambda r_k}$ for every $\lambda \in \mathbb{N}$.
By Remark \ref{P''}.(1), $P_{n+\lambda r_k}'=P_n'$ and $P_{n+\lambda r_k}''=\{i+\lambda d r_k \mid i \in P_n''\}$ for every $\lambda \in \mathbb{N}$. Therefore, by Theorem \ref{PF}, $\psi_n^\lambda(i)=i$ if $i \in P_n'$, whereas $\psi_n^\lambda(i)=i+\lambda dr_k$ if $i \in P_n''$. It immediately follows that $N_1$ and $N_3$ are the same for $M_n$ and $M_{n+\lambda r_k}$. As for $N_2$, it is enough to note that Remark \ref{OP}.(1) implies $\m(i+\lambda dr_k)=\m(i)+\lambda d$ if $i \in P_{n}''$ because $i\geq dn-r_k >r_k^2-r_k$. 
\end{remark}

In the following remark, we show that $N$ is always less than $r_k^4$. Therefore, for simplicity, in the statements of the results in this section, one can replace $N$ with $r_k^4$. However, in general, $N$ provides a much better bound, as shown in Example \ref{example NG}.

\begin{remark} \rm
The inequality $N < r_k^4$ always holds. To prove this, we assume that $n>N_0$.  As noted in the proof of Proposition \ref{PF ordered},   $N_0<r_k^4$. We divide the proof into four parts.
\begin{itemize}[leftmargin=*]
\item Let $i,\ell \in P_n'$. Then, by Remark \ref{P''}.(3), $r_k+\ell-i<r_k+d(r_k^3-r_k^2)<dr_k^3<r_k^4$. Moreover,
\[
(\m(\ell)-\m(i))r_k-\ell+i\leq \left(\frac{\ell}{r_1}-\frac{i}{r_k}\right)r_k+i=\ell \frac{r_k}{r_1} \leq \ell \frac{r_k}{d} < d(r_k^3-r_k^2)\frac{r_k}{d}<r_k^4.
\]
It follows that $N_1 < r_k^4$.
\item Let $i,\ell \in P_n''$. We have
\[
(\m(\ell)-\m(i))r_k-\ell+i\leq \left(\m(\ell)-\frac{i}{r_k}\right)r_k-\ell+i=\m(\ell)r_k-\ell.
\]
Let $\lambda$ be the maximum integer such that $dn-\lambda r_k >N_0$.
In particular, using the maximality of $\lambda$ together with (3) and (4) of Remark \ref{OP}, we obtain 
$$\ell -(\lambda+1)r_k \leq dn+\F(S)-(\lambda+1)r_k \leq N_0+\F(S).$$
If $\F(S)<0$, then $N_0+\F(S)=r_k^2+\F(S)< 2r_k^2-2r_k\leq r_k^3-2r_k$.
Also, when $\F(S)>0$, we have $N_0+\F(S)=r_k^2+(r_k+1)\F(S)< r_k^2+(r_k+1)(r_k^2-2r_k)=r_k^3-2r_k$
 Hence, either way, $\ell -(\lambda+1)r_k  < r_k^3-r_k$. 

Moreover, $\ell-\lambda r_k \geq dn-r_k-\lambda r_k>N_0-r_k=r_k^2-r_k+\max\{0, \F(S)r_k\}  \geq r_k^2-r_k$.
Then, by Remark \ref{OP}.(1), we have
\begin{gather*}
\m(\ell)r_k-\ell=\m(\ell-\lambda r_k)r_k-(\ell-\lambda r_k)< (\ell-\lambda r_k) r_k= \\
= (\ell-(\lambda+1) r_k) r_k+r_k^2< (r_k^3-r_k)r_k+r_k^2=r_k^4.
\end{gather*}
Therefore, $N_2<r_k^4$.
\item Let $\ell \in P_n'$. Again by Remark \ref{P''}.(3), since $r_k \geq 2$, we get
\[
\m(\ell)r_k+3r_k \leq \frac{\ell}{r_1}r_k+3r_k < d(r_k^3-r_k^2)\frac{r_k}{d}+3r_k=r_k^4-r_k^3+3r_k < r_k^4. 
\]
In particular, it follows that $N_3<r_k^4$.
\item Finally, by Remark \ref{OP}.(4), we have $N_0+r_k+\F(S)/d < r_k^2+ r_k(r_k^2-2r_k)+  r_k+(r_k^2-2r_k)
=r_k^3-r_k<r_k^4$.
\end{itemize}
By putting together the four parts, we obtain $N < r_k^4$.
\end{remark}

Given a numerical semigroup $H$, the {\it canonical ideal} of $H$ is defined as $K(H)=\{x \in \mathbb{N} \mid \F(H)-x \notin H\}$, while the {\it trace ideal} of $K(H)$ is $\tr(H)=K(H)+(H-K(H))$, where $H-K(H)=\{x \in \mathbb{Z} \mid x+K(H) \subseteq H\}$. It follows from the definition that $\tr(H) \subseteq H$, and it is well known that $\tr(H)=H$ if and only if $H$ is symmetric, i.e., $H=K(H)$. In \cite{HHS, HHS1} the authors introduce a new notion that generalizes symmetric numerical semigroups: $H$ is said to be {\it nearly Gorenstein} if $\tr(H)$ contains the maximal ideal $H\setminus \{0\}$. This is equivalent to saying that the associated numerical semigroup ring $\Bbbk[[H]]$ is nearly Gorenstein as defined in \cite{HHS}.  

Rather than the definition above, in this paper we will use a characterization established in \cite{MS}. Let $H=\langle h_0, h_1, \dots, h_k \rangle$ be a numerical semigroup minimally generated by $h_0, h_1, \dots, h_k$. We say that $(f_0, f_1, \dots, f_k) \in \PF(H)^{k+1}$ is a {\it nearly Gorenstein vector} for $H$,  an NG-vector for short, if $h_i+f_i-f \in H$ for every $f \in \PF(H)$ and every $i=0,1, \dots, k$. By \cite[Proposition 1.1]{MS}, a numerical semigroup $H$ is nearly Gorenstein if and only if there exists an NG-vector for $H$. Using this notion, in the following theorem we will prove that nearly Gorensteinness is eventually a periodic property in the shifted family $M_n$.

\begin{theorem} \label{NearlyGorenstein} 
If $n > N$ and $M_n$ is a nearly Gorenstein numerical semigroup, then $M_{n+\lambda r_k}$ is nearly Gorenstein for every $\lambda \in \mathbb{N}$. Moreover, if $(f_0, \dots, f_{k})$ is an {\rm NG}-vector for $M_n$, then $(\varphi_n^\lambda(f_0), \dots, \varphi_n^\lambda(f_{k}))$ is an {\rm NG}-vector for $M_{n+\lambda r_k}$ for all $\lambda \in \mathbb{N}$.
\end{theorem}

\begin{proof}
Assume that $M_n$ is nearly Gorenstein and has NG-vector $(f_0, \dots, f_{k})$. In light of Remark \ref{N}, it is enough to prove that $M_{n+r_k}$ is nearly Gorenstein with NG-vector $(\varphi_n(f_0), \dots, \varphi_n(f_{k}))$.

By definition, $n+r_j+f_j-f \in M_n$ for every $f \in \PF(M_n)$ and $j=0,1,\dots,k$. It is enough to prove that $n+r_j+r_k+\varphi_n(f_j)-\varphi_n(f) \in M_{n+r_k}$ for every $f \in \PF(M_n)$ and every $j$; indeed, this means that $(\varphi_n(f_0), \dots, \varphi_n(f_{k}))$ is an NG-vector for $M_{n+r_k}$.

Fix $f\in \PF(M_n)$ and $j \in \{0,1,\dots, k\}$. In particular, $f+n \in \Ap(M_n,n)$, then by Theorem \ref{Apéry} we may write $f=i+(\m(i)-1)n$ and $f_j=\ell+(\m(\ell)-1)n$ for some $i,\ell \in \Ap(S,dn)$. We set 
\begin{gather*}
\mathcal A:=n+r_j+f_j-f=n+r_j+\ell-i+(\m(\ell)-\m(i))n \\
\mathcal B:=n+r_j+r_k+\varphi_n(f_j)-\varphi_n(f).
\end{gather*}
By assumption $\mathcal A \in M_n$ and we only need to prove that $\mathcal B \in M_{n+r_k}$.
We distinguish four cases.

\medskip
 
{\bf Case 1: $i,\ell < dn-r_k$.} By Corollary \ref{CorPF} we have
\[
\mathcal B=n+r_j+f_j -f + r_k + (\m(\ell)-\m(i)) r_k=\mathcal A+(\m(\ell)-\m(i)+1) r_k.
\]
If we prove that $\mathcal A$ has a factorization with length $\m(\ell)-\m(i)+1$ in $M_n$, then by adding $r_k$ to every generator we get a factorization of $\mathcal B$ with length $\m(\ell)-\m(i)+1$ in $M_{n+r_k}$, and in particular $\mathcal B \in M_{n+r_k}$. 

By $r_j+\ell-i\leq N_1 \leq N<n$, we get $\mathcal{A}<(\m(\ell)-\m(i)+2)n$. Since $n$ is the smallest generator of $M_n$, we obtain that every factorization of $\mathcal{A}$ in $M_n$ has length at most $\m(\ell)-\m(i)+1$.

Moreover, $n>N_1 \geq (\m(\ell)-\m(i))r_k-\ell+i$, then we have $n+r_j+\ell-i>(\m(\ell)-\m(i))r_k$. It follows that $\mathcal{A}>(\m(\ell)-\m(i))(n+r_k)$, which means that every factorization of $\mathcal{A}$ in $M_n$ has length at least $\m(\ell)-\m(i)+1$ because $n+r_k$ is the biggest generator of $M_n$. Since $\mathcal{A} \in M_n$, it has a factorization and its length has to be $\m(\ell)-\m(i)+1$, as we needed to prove.

\medskip

{\bf Case 2: $i,\ell \geq dn-r_k$.} As in the previous case, we have
\[
\mathcal B=n+r_j+f_j -f + r_k + (\m(\ell)-\m(i)) r_k=\mathcal A+(\m(\ell)-\m(i)+1) r_k
\]
and we only need to find a factorization of $\mathcal{A}$ in $M_n$ having length $\m(\ell)-\m(i)+1$. By using (3) and (4) of Remark \ref{OP}, we get 
\[
r_j+\ell-i\leq r_j+dn+\F(S)-(dn-r_k)=r_j+\F(S)+r_k \leq 2r_k+\F(S) \leq r_k^2\leq N_0 <N<n.
\] 
Therefore, $\mathcal{A}<(\m(\ell)-\m(i)+2)n$ and so every factorization of $\mathcal{A}$ in $M_n$ has length at most $\m(\ell)-\m(i)+1$.

Since $n>N_2 \geq (\m(\ell)-\m(i))r_k-\ell+i$, it follows that $n+r_j+\ell-i>(\m(\ell)-\m(i))r_k$ and, as in the previous case, $\mathcal{A}>(\m(\ell)-\m(i))(n+r_k)$. This implies that every factorization of $\mathcal{A}$ in $M_n$ has length exactly $\m(\ell)-\m(i)+1$. 

\medskip

{\bf Case 3: $i \geq dn-r_k$ and $\ell < dn-r_k$.} We claim that this case cannot occur. Since $n>N_3 \geq (\m(\ell)r_k+3r_k)/d$, we get
\[
\m(i) \geq \frac{i}{r_k} \geq \frac{dn-r_k}{r_k} > \frac{\m(\ell)r_k+2r_k}{r_k}=\m(\ell)+2.
\]
This means that $\m(\ell)-\m(i)+2 < 0$. On the other hand, since $\ell-i<0$, it follows that $\mathcal{A}<(\m(\ell)-\m(i)+2)n<0$ and this is not possible because $\mathcal{A} \in M_n$. 

\medskip

{\bf Case 4: $i < dn-r_k$ and $\ell \geq dn-r_k$.} By using $n>N \geq N_0+r_k+\F(S)/d$ and Remark \ref{P''}.(2), we have
\[r_j+\ell-i \geq r_j+dn-r_k-i> 
r_j+d(N_0+r_k)+\F(S)-r_k-(dN_0+dr_k-r_k) =r_j+\F(S)\geq \F(S),\]
thus $r_j+\ell-i \in S$. If $r_j+\ell-i \in \Ap(S,dn)$, then by Theorem \ref{Apéry} $r_j+\ell-i + \m(r_j+\ell-i)n \in \Ap(M_n,n)$ and it has a factorization in $M_n$ with length $\m(r_j+\ell-i)$; this implies that $\m(r_j+\ell-i) \leq \m(\ell)-\m(i)+1$ because $\mathcal{A} \in M_n$. Moreover, since
\[
\mathcal{A}=r_j+\ell-i + \m(r_j+\ell-i)n+(\m(\ell)-\m(i)+1-\m(r_j+\ell-i))n
\]
and $n$ is a generator of $M_n$, we get a factorization of $\mathcal A$ in $M_n$ with length $\m(\ell)-\m(i)+1$.

On the other hand, if $r_j+\ell-i \notin \Ap(S,dn)$, then it is larger than $dn$ and using (3) and (4) of Remark \ref{OP} we obtain 
\[0 \leq r_j+\ell-i-dn \leq r_j+dn+\F(S)-dn<r_j+r_k^2-2r_k<r_k^2<dn.\] 
Hence, $r_j+\ell-i-dn \in \Ap(S,dn)$ and $r_j+\ell-i-dn + \m(r_j+\ell-i-dn)n \in \Ap(M_n,n)$ has a factorization with length $\m(r_j+\ell-i-dn)$. Therefore,
\[
\mathcal{A}=r_j+\ell-i-dn + \m(r_j+\ell-i-dn)n+(\m(\ell)-\m(i)+d+1-\m(r_j+\ell-i-dn))n
\]
yields a factorization of $\mathcal{A}$ with length $\m(\ell)-\m(i)+d+1$.

In any case, there is a factorization of $\mathcal{A}$ in $M_n$ with length either $\m(\ell)-\m(i)+1$ or $\m(\ell)-\m(i)+d+1$. Moreover, by Theorem \ref{PF}
\begin{align*}
\mathcal B &=n+r_j+f_j -f + r_k + (\m(\ell)-\m(i)+2d) r_k+dn= \\ 
&=\mathcal A+(\m(\ell)-\m(i)+1) r_k + d(n+2r_k)= \\
&=\mathcal A+(\m(\ell)-\m(i)+d+1) r_k + d(n+r_k).
\end{align*}
Since both $n+r_k$ and $n+2r_k$ are generators of $M_{n+r_k}$, in both situations we get $\mathcal{B} \in M_{n+r_k}$.
\end{proof}

\begin{example} \label{example NG} \rm
Consider the shifted family $M_n=\langle n, n+2, n+3, n+5\rangle$. For $M_{40}=\langle 40,42,43,45 \rangle$ we find with GAP that $\PF(M_{40})=\{359, 361\}$. It is easy to check that $M_{40}$ is nearly Gorenstein with NG-vector $(361,359,361,359)$. Moreover, $r_1=2, r_2=3, r_3=5, d=1, S=\langle 2,3 \rangle$ and
\begin{align*}
&359=39+9 \cdot 40 - 40 &\text{where } 39 \in \Ap(S,40) \text{ and } \m(39)=9,\\
&361=41+9 \cdot 40 - 40 &\text{where } 41 \in \Ap(S,40) \text{ and } \m(41)=9.
\end{align*}
Since $39, 41 \in P_{40}''$, we have $N_1=N_3=0$, $N_2=2$ and then $N=r_3^2+(\F(S)+1)r_3+\F(S)=25+10+1=36$. Note that $N$ is much less than $r_3^4=625$. Since $40>N$, Theorem \ref{NearlyGorenstein} and Corollary \ref{lambda} imply that $M_{40+5\lambda}$ is a nearly Gorenstein numerical semigroup with NG-vector 
\[(361+85\lambda+5\lambda^2,359+85\lambda+5\lambda^2,361+85\lambda+5\lambda^2,359+85\lambda+5\lambda^2),\]
for all $\lambda \in \mathbb{N}$. 
\end{example}

\section{Almost symmetric property}

Almost symmetric numerical semigroups were introduced by Barucci and Fr\"oberg in \cite{BF}, together with the more general notion of almost Gorenstein rings, and have been much studied in the last decades. They are defined as the numerical semigroups $H$ for which $M(H) + K(H) \subseteq M(H)$, where $M(H)=H \setminus \{0\}$ is the maximal ideal of $H$. It turns out that almost symmetric numerical semigroups are always nearly Gorenstein and, more precisely, they are exactly the nearly Gorenstein numerical semigroups that admit $(\F(H), \F(H), \dots, \F(H))$ as NG-vector, see \cite[Proposition 1.1]{MS}, \cite[Proposition 6.1]{HHS}. Therefore, Proposition \ref{Frobenius} and Theorem \ref{NearlyGorenstein} immediately imply the following result.

\begin{corollary} \label{Almost symmetric}
Let $n > N$. If $M_n$ is almost symmetric, then $M_{n+\lambda r_k}$ is almost symmetric for every $\lambda \in \mathbb{N}$.
\end{corollary}

Another useful characterization of almost symmetric numerical semigroups was given by Nari in \cite[Theorem 2.4]{N}, who showed that if $f_1<f_2< \dots < f_t=\F(H)$ are the pseudo-Frobenius numbers of a numerical semigroup $H$, then $H$ is almost symmetric if and only if $f_a+f_{t-a}=\F(H)$ for every $a=1, \dots, t-1$. 

\begin{example} \rm
In \cite[Example 7.7]{HW}, Herzog and Watanabe stated without proof that the semigroup $\langle 10+4 \lambda, 11+4\lambda,13+4\lambda,14+4\lambda \rangle$ is almost symmetric of type $3$ for every $\lambda \in \mathbb{N}$. They also wrote that this example is due to T. Numata. To see why this is true, we first consider the shifted family $M_n=\langle n, n+1, n+3, n+4\rangle$, where $r_1=1, r_2=3, r_3=4, d=1$, $S=\langle 1,3, 4\rangle=\N,\F(S)=-1$.
Then $M_{10+4\lambda}=\langle 10+4 \lambda, 11+4\lambda,13+4\lambda,14+4\lambda \rangle$.

 Although Corollary \ref{lambda} does not apply directly to $M_{10}$ since $N_0=16$, nevertheless, it can be applied to $M_{18}$.
Observe that $\PF(M_{18})=\{20,69,89\}$. One can readily verify that the resulting formulas are also valid for $M_{10}$ and $M_{14}$ as well. In particular, we obtain 
\[\PF(M_{10+4\lambda})=\PF(M_{18+4(\lambda-2)})=
\{f_1=12+4\lambda, f_2=17+18\lambda+4\lambda^2,f_3=29+22\lambda+4\lambda^2\}
\] 
for every $\lambda \in \mathbb{N}$. Since it always holds that $f_1+f_2=f_3$, it follows that $M_{10+4\lambda}$ is almost symmetric of type $3$ for every $\lambda \in \mathbb{N}$.
\end{example}

Despite the preceding example, Herzog and Watanabe proved in \cite[Theorem 7.6]{HW} that if $M_n$ has embedding dimension $4$ and $n$ is large enough, then $M_{n}$ is not almost symmetric of type $2$. In the following proposition, we generalize this result by considering arbitrary embedding dimension and all even types.

\begin{proposition} \label{even-type}
If $n \geq r_k^4$, then $M_{n}$ is not almost symmetric of even type. 
\end{proposition}

\begin{proof} Assume by contradiction that $M_n$ is almost symmetric of even type $t$. 
Let $f_1<f_2< \dots < f_t$ be the pseudo-Frobenius numbers of $M_n$. Since $M_n$ is almost symmetric, by Nari's characterization \cite[Theorem 2.4]{N}, $f_a+f_{t-a}=\F(M_n)$ for every $a=1, \dots, t-1$, and in particular $2f_{t/2}=\F(M_n)$ because $t$ is even. 

By Remark \ref{N} and Corollary \ref{Almost symmetric}, also $M_{n+ \lambda r_k}$ is almost symmetric of type $t$ for every $\lambda \in \mathbb{N}$. Furthermore, Proposition \ref{PF ordered} together with Nari's characterization imply $2\varphi_n^\lambda(f_{t/2})=\varphi_n^\lambda(\F(M_n))$.
Let $\F(M_n)=i+\m(i)n-n$ and $f_{t/2}=j+\m(j)n-n$ with $i,j \in \Ap(S,dn)$. By Proposition \ref{Frobenius}, $\varphi_n^\lambda(\F(M_n))=
\F(M_n)+(\m(i)+(\lambda +1)d-1) \lambda r_k+\lambda dn$.

If $j <dn-r_k$, by Corollary \ref{lambda}, we have
\[
2f_{t/2}+2(\m(j)-1)\lambda r_k=\F(M_n)+(\m(i)+(\lambda+1)d-1) \lambda r_k+\lambda dn,\]
which holds only if either $\lambda=0$ or
\[
\lambda=\frac{(2\m(j)-\m(i)-d-1) r_k-dn}{d r_k}.
\]
Since the equality should be true for every $\lambda \in \mathbb{N}$, we get a contradiction.

If $j \geq dn-r_k$, again by Corollary \ref{lambda}, we have
\[
2f_{t/2}+2(\m(j)+(\lambda+1)d-1) \lambda r_k+2\lambda dn=\F(M_n)+(\m(i)+(\lambda+1)d-1) \lambda r_k+\lambda dn,
\]
which holds when either $\lambda=0$ or
\[
\lambda=\frac{(\m(i)-2\m(j)-d+1) r_k- dn}{d r_k}.
\]
Thus, we obtain a contradiction when $\lambda$ differs from these values.
\end{proof}

Notice that the previous result is false if we replace almost symmetric with nearly Gorenstein. See, for instance, Example \ref{example NG}.

\section{Residue}

A notion closely related to the nearly Gorenstein property is the {\it residue} of a numerical semigroup $H$, defined as $\res(H)=|H\setminus \tr(H)|$. Clearly, $\res(H)=0$ if and only if $H$ is symmetric, while $\res(H) \leq 1$ if and only if $H$ is nearly Gorenstein; then, the residue is a measure of how far a numerical semigroup is from being symmetric. Since symmetric and nearly Gorenstein properties are eventually periodic in a shifted family, one may expect that also the residue is periodic. This was verified for $k=2$ in \cite[Theorem 4.2]{HHS1}, but we will see in Example \ref{res} that this is not usually the case, even if the type of the semigroup is two. Before illustrating this example, we need to state a proposition that we will use. This is a special case of a more general result proved in \cite[Proposition 1.5]{JSZ} for simplicial affine semigroups, but we include a simple proof here for the sake of completeness.

\begin{proposition}\label{JSZ} For every numerical semigroup $H$ it holds
$$\tr(H)=\{h \in H \mid \exists \, f^* \in \PF(H) \text{ such that } h+f^*-f \in H \ \forall f \in \PF(H)\}.$$
\end{proposition}

\begin{proof}
We recall that the relative ideal $K(H)$ is generated by the elements $\F(H)-f$ with $f \in \PF(H)$. Assume first that $h \in \tr(H)=K(H)+(H-K(H))$, so that $h=x+y$ with $x \in K(H)$ and $y \in (H-K(H))$. Thus, $x=\F(H)-f^*+z$ for some $f^* \in \PF(H)$ and $z \in H$. It follows that $h+f^*-\F(H)=y+z \in (H-K(H))$, and then $h+f^*-f=(h+f^*-\F(H))+(\F(H)-f) \in H$ for every $f \in \PF(H)$.

Conversely, write $h$ as $h=(\F(H)-f^*)+(h+f^*-\F(H))$. Note that $\F(H)-f^* \in K(H)$, thus it is enough to prove that $(h+f^*-\F(H)) \in (H-K(H))$. To verify this, observe that for every $f \in \PF(H)$, we have $(h+f^*-\F(H))+(\F(H)-f)=h+f^*-f \in H$ by hypothesis, which immediately implies $h+f^*-\F(H) \in (H-K(H))$.
\end{proof}

\begin{example}\label{res} \rm
Consider the shifted family $M_n=\langle n, n+2, n+3, n+7 \rangle$, where $r_1=2, r_2=3, r_3=7, d=1, S=\langle 2,3\rangle$  has $\F(S)=1$ and $N_0=56$.  

For brevity, we set $p=63+7\lambda, m=p+7=70+7\lambda$, where $\lambda \in \N$.
We will show that
$$M_{p}\setminus \tr(M_p)=\{0,m,2m \dots, (\lambda+8)m\};$$
in particular, $\res(M_{63+7\lambda})=\lambda+9$ and so the residue of $M_{63+7\lambda}$ grows linearly with $\lambda$. Recalling that $M_p=\langle p,p+2, p+3, m=p+7 \rangle$ and that $\tr(M_p)$ is an ideal of $M_p$, in order to prove the equality above we only need to show that
\begin{gather}\label{elements for residue}
p, p+2, p+3, (\lambda +9)m \in \tr(M_p) \hspace{20pt} {\rm and} \hspace{20pt} (\lambda+8)m \notin \tr(M_p). 
\end{gather}  
We will make use of Proposition \ref{JSZ}, so we first note that $\PF(M_{63})=\{627,694\}$ with $63>N_0$ and, by Corollary \ref{lambda}, it easily follows that 
$$\PF(M_{p})=\{f_1=627+133 \lambda+7\lambda^2, f_2=694+140\lambda+7\lambda^2\}.$$ 
Note that the difference between the two pseudo-Frobenius numbers of $M_{p}$ is $f_2-f_1=67+7\lambda=p+4$, and therefore an element $x \in M_{p}$ is in $\tr(M_{p})$ if and only if either $x+p+4$ or $x-(p+4)$ is in $M_{p}$. Now we analyze each element in (\ref{elements for residue}):
\begin{itemize}[leftmargin=*]
\item $p\in \tr(M_p)$ because $p+p+4=2(p+2) \in M_p$.
\item $p+2 \in \tr(M_p)$ because $p+2+p+4=2(p+3) \in M_p$.
\item $p+3\in \tr(M_p)$ because $p+3+p+4=p+(p+7) \in M_p$.
\item $(\lambda +9)m \in \tr(M_p)$ because $(\lambda +9)m+p+4=(\lambda +9)(p+7)+p+4=(\lambda +9)p+2(p+2) \in M_p$.
\item It remains to show that $x=(\lambda+8)m \notin \tr(M_p)$, that is $x+p+4,x-(p+4) \notin M_p$. These hold because 
\begin{align*}
&x+p+4=(\lambda +8)(70+7\lambda)+63+7\lambda+4=627+133\lambda+7\lambda^2=f_1 \notin M_p, \\
&x-(p+4)=f_1-2(p+4)=f_2-3(p+4)=f_2-((p+2)+(p+3)+(p+7)) \notin M_p,
\end{align*}
where the last element is not in $M_p$ because $(p+2)+(p+3)+(p+7) \in M_p$ and $f_2 \notin M_p$.
\end{itemize}
\end{example}

\begin{example}\label{res-2} \rm
With some computation, one can also see that $M_{46+11\lambda}=\langle 46+11\lambda, 48+11\lambda,52+11\lambda,57+11\lambda \rangle$ has type $3$ and residue $\lambda +8$. Indeed, if we set $p=46+11\lambda$, the elements in the semigroup $M_p$ that are not in the trace ideal are exactly $0,p,p+11,p+(p+11),2p,2p+(p+11),2(p+11), \dots, (\lambda+3)(p+11)$.
\end{example}

Nevertheless, in all examples we know, the residue appears to exhibit linear behaviour. This naturally leads to the following question.

\begin{question} \rm
Given a numerical semigroup $M_n$, is $\res(M_{n+\lambda r_k})$ linear in $\lambda$ for $\lambda \gg 0$? 
\end{question}

\section{More results}

\subsection{Canonical reductions}

A one-dimensional Cohen-Macaulay local ring is said to have a {\it canonical reduction} if it admits a canonical ideal that is a reduction of the maximal ideal, see \cite{R}. We say that a numerical semigroup $H$ has a canonical reduction if the semigroup ring $\Bbbk[[H]]$ does. From \cite[Theorem 3.13]{R} it immediately follows that $M_n$ has a canonical reduction if and only if $n+\F(M_n)-f \in M_n$ for every $f \in \PF(M_n)$. These semigroups were also introduced independently in \cite{BFR}, where they are called {\it positioned numerical semigroups}.

\begin{corollary} \label{canonical reduction}
If $n > N$ and $M_n$ has a canonical reduction, then $M_{n+\lambda r_k}$ has a canonical reduction for every $\lambda \in \mathbb{N}$.
\end{corollary}

\begin{proof} 
Since $\F(M_{n+r_k})=\varphi_n(\F(M_n))$ by Proposition \ref{Frobenius}, this is a consequence of the proof of Theorem \ref{NearlyGorenstein} when one chooses $j=0$ and $f_0=\F(M_n)$.  
\end{proof}

An explicit example is given below.

\begin{example} \rm
We consider the shifted family $M_n=\langle n, n+2, n+3, n+4 \rangle$, and in particular $M_{26}$. In this case, we have $\PF(M_{26})=\{181,183\}$ and $P_{26}=P''_{26}=\{25,27\}$, so that $N_1=N_3=0$. Moreover, $\m(25)=\m(27)=7$ implies $N_2=2$, and hence $N=25$. Note that $M_{26}$ has a canonical reduction: indeed, $\F(M_{26})=183$ and there is only one other pseudo-Frobenius number, thus we only need to check that $26+183-181=28 \in M_{26}$, which is true. Since $26>N$, Corollary \ref{canonical reduction} implies that $M_{26+4\lambda}$ has a canonical reduction for every $\lambda \in \mathbb{N}$. Note also that $M_{26}$ is not nearly Gorenstein, because $29+183-181=31 \notin M_{26}$ and $29+181-183=27 \notin M_{26}$.
\end{example}

\subsection{Reduced type}

Let $R$ be a non-regular one-dimensional complete local domain with residue field $\Bbbk$, which is also a $\Bbbk$-algebra.
Huneke, Maitra, and Mukundan introduced a new invariant, the {\it reduced type} $s(R)$, in order to study Berger's conjecture; see \cite{HMM} for more details. In \cite{MM}, the last two authors studied the reduced type in the specific case of numerical semigroup rings, and proved that, for a numerical semigroup $H$, the reduced type of $\Bbbk[[H]]$ is equal to the number of pseudo-Frobenius numbers of $H$ in the interval $[\F(H)-n,\F(H)]$, where $n$ is the smallest nonzero integer in $H$, see \cite[Theorem 2.13]{MM}. In the following proposition, we show that also the reduced type is eventually periodic in a shifted family of numerical semigroups.

\begin{proposition} \label{r-type}
Assume $n \geq r_k^4$. Then, $s(\Bbbk[[M_{n}]])=s(\Bbbk[[M_{n+\lambda r_k}]])$ for every $\lambda \in \mathbb{N}$. Moreover, if $\F(M_n)=i+\m(i)n-n$, with $i \in \Ap(S,dn)$, then
\[
s(\Bbbk[[M_{n}]])=|\{j \in P_n'' \mid \m(j)=\m(i)\}| + |\{j \in P_n'' \mid \m(j)=\m(i)-1 \text{ and } j>i \}| \leq |P_n''|.
\]
\end{proposition}

\begin{proof}
Let $\F(M_n)=i+\m(i)n-n$ and let $f=j+\m(j)n-n$ be another pseudo-Frobenius number of $M_n$, where $i,j \in \Ap(S,dn)$. By Proposition \ref{Frobenius}, we have $i \in P_n''$ and $\F(M_{n+\lambda r_k})=\varphi_n^\lambda(\F(M_n))$. If $j \in P_n'$, in the proof of Proposition \ref{PF ordered} we proved that $\m(i)>\m(j)$ and then
\[
\F(M_n)-f=i-j+(\m(i)-\m(j))n>n.
\]
Bearing in mind Remark \ref{P''}, this means that we can focus only on the case $j \in P_n''$, since the pseudo-Frobenius numbers corresponding to elements of $P_n'$ do not give any contribution for the reduced type. 

When $j \in P_n''$, in the proof of Proposition \ref{PF ordered} we have seen that $\m(i) \geq \m(j)$. We distinguish three cases, and in all of them we use Remark \ref{OP} and Corollary \ref{lambda}.
\begin{itemize}[leftmargin=*]
\item $\m(j)=\m(i)$. Then, $\varphi_n^\lambda(\F(M_n))-\varphi_n^\lambda(f)=\F(M_n)-f=i-j<(dn+\F(S))-(dn-r_k)=\F(S)+r_k<n$. 
\item $\m(j)=\m(i)-1$. Then, $\varphi_n^\lambda(\F(M_n))-\varphi_n^\lambda(f)=\F(M_n)-f+\lambda r_k=i-j+n+\lambda r_k$. This means that $\varphi_n^\lambda(f)$ is in the interval $[\F(M_{n+\lambda r_k})-(n+\lambda r_k),\F(M_{n+\lambda r_k})]$ if and only if $f$ is in the interval $[\F(M_n)-n,\F(M_n)]$, and this happens precisely when $j>i$.
\item $\m(j)\leq \m(i)-2$. Then, we have
\begin{gather*}
\F(M_n)-f\geq i-j+2n \geq (dn-r_k)-(dn+\F(S))+2n=2n-(r_k+\F(S))>n \\
\varphi_n^\lambda(\F(M_n))-\varphi_n^\lambda(f)\geq \F(M_n)-f+2\lambda r_k >n+\lambda r_k,
\end{gather*}
thus in this case we never have contributions for the reduced type.
\end{itemize} 
Putting everything together, we obtain the desired result. 
\end{proof}

\begin{remark} \rm The formula in Proposition \ref{r-type} is proved to hold for $n\geq r_k^4$. Let $\mu$ be a positive integer such that $p:=n-\mu r_k > N_0$. Then  $P_n^{''}=\{j+\mu d r_k \mid j \in P_{n-\mu r_k}^{''} \}$. Using  Corollary \ref{shift-m}, Propositions \ref{r-type} and \ref{PF ordered}, we get that
$$
s(\Bbbk[[M_{n}]])= |\{j \in P_p{''} \mid \m(j)=\m(i)\}| + |\{j \in P_p'' \mid \m(j)=\m(i)-1 \text{ and } j>i \}| \leq |P_p''|.
$$
\end{remark}

We illustrate with an example.

\begin{example}  \rm
We return to the shifted family $M_n=\langle n, n+2, n+7, n+11\rangle$ from Example \ref{Counterexample}, where $N_0=176$ and $r_3^4=11^4=14641$.
To find the reduced type of $\Bbbk[[M_{14643}]]$, since $14643 \equiv 200 \mod 11$,  it is enough to inspect $P_{200}^{''}=\{191, 205, 195, 199, 203\}$. Its elements have minimal decomposition length $19, 19, 20, 20, 20$, respectively. As $\F(M_{200})=4003 \equiv 203\mod 200$, we get that $s(\Bbbk[[M_{14643}]])=|\{195, 199, 203\}|+|\{205\}|=4$.
\end{example}

\end{document}